\documentclass[11pt,leqno, openbib]{article}

\usepackage{amsthm}
\usepackage{amsmath}
\usepackage{amssymb}
\usepackage{amsfonts}
\usepackage{amscd}
\usepackage{palatino}
\usepackage[mathscr]{eucal}
\usepackage{epsfig}
\usepackage{verbatim}
\usepackage{latexsym}
\usepackage{graphics}

\begin{document}

\newcounter{chn}
\newenvironment{nch}[1]
               {\addtocounter{chn}{1}
               \newpage
                      \begin{center}
                      {\Large {\bf {#1}}}
                      \end{center}
                      \vspace{5mm}}{}

%The Theorem env.
\newcounter{thn}
\renewenvironment{th}
{\addtocounter{thn}{1}
 \em 
     \noindent 
     {\sc Theorem \arabic{thn} . ---}}{}

%Theorem without numbers
\newenvironment{theo}
   {\em
        \vspace{2mm}
        \noindent
                {\sc Theorem\/. ---}}{}

%Lemma without numbers
\newenvironment{lem}
   {\em
        \vspace{2mm}
        \noindent
                {\sc Lemma\/. ---}}{}

%Prop without numbers
\newenvironment{pro}
   {\em
        \vspace{2mm}
        \noindent
                {\sc Proposition\/. ---}}{}

%The Lemma env.
\newcounter{len}
\renewenvironment{le}
{\addtocounter{len}{1}
 \em 
     \vspace{2mm}
     \noindent 
     {\sc Lemma \arabic{len}. ---}}{}

%The Proposition env.
\newcounter{prn}
\newenvironment{prop}
{\addtocounter{prn}{1}
  \em 
      \vspace{2mm}
      \noindent 
      {\sc Proposition \arabic{prn}. ---}}{}

\newcounter{exn}
\newenvironment{prob}
                     {\addtocounter{exn}{1} 
                      \vspace{2mm}
                      \noindent 
                      {\em Exercise \arabic{exn} :}}{}

\newcounter{corn}
\newenvironment{cor}
  {\em 
     \addtocounter{corn}{1}
     \vspace{2mm}   
     \noindent 
        {\sc Corollary \arabic{corn}. ---}}{\vspace{1mm}}

%Corollary without numbers
\newenvironment{cor1}
   {\em 
        \vspace{2mm}
        \noindent
                {\sc Corollary\/. ---}}{}

%Definitions
\newcounter{defn}
\newenvironment{df}{
                    \vspace{2mm}
                    \addtocounter{defn}{1} 
                    \noindent 
                            {\sc Definition \arabic{defn}. ---}}{\vspace{1mm} }{}
%paragraphs
\newcounter{pgn}
\newenvironment{pg}{
                    \vspace{2mm}
                    \addtocounter{pgn}{1}
                    \noindent
                            {\arabic{scn}.\arabic{pgn}\;}}{}

%The section env.
\newcounter{scn}
\newenvironment{nsc}[1]
              {\addtocounter{scn}{1}
              \vspace{3mm}
                     \begin{center}
                     {\sc { \arabic{scn}. {#1}}}
                     \end{center}
              \setcounter{equation}{0}
              \setcounter{prn}{0}
              \setcounter{len}{0}
              \setcounter{thn}{0}
              \setcounter{pgn}{0}
              \setcounter{corn}{0}
              \vspace{3mm}}{}

\begin{center}
\Large {\bf 

{
Arithmetical Applications of an Identity for the  Vandermonde Determinant. 
}
}

\end{center}
\vspace{5mm}

\vspace{-10mm}
\begin{center}
\it{

by \\

\vspace{3mm}

D.S. Ramana

}
\end{center}

\begin{center}
{\bf Abstract}
\end{center}

\vspace{1mm}
\noindent
\vspace{1mm}
\noindent
When $\{\alpha_i\}_{1 \leq i \leq m}$ is a sequence of distinct non-zero elements of an integral domain $A$ and $\gamma$ is a common multiple of the $\alpha_i$ in $A$ we obtain, by means of a simple identity for the
 Vandermonde determinant, a lower bound for $\sup_{1\leq i < j \leq m}\phi(\alpha_i - \alpha_j)$ in terms of $\phi(\gamma)$, where $\phi$ is a function from the nonzero elements of $A$ to ${\bf R}_{+}$ satisfying certain natural conditions. We describe several applications of this bound.

\vspace{2mm}
\begin{nsc}{Introduction}

\noindent
This article is concerned with the following question. Suppose that $\{\alpha_i\}_{1 \leq i \leq m}$ is a sequence of distinct elements in an integral domain $A$ and that all the $\alpha_i$ have a common multiple $\gamma \neq 0 $ in $A$. Let $\phi$ be a function from $A$ into ${\bf R}_{+}$ satisfying $\phi(xy) = \phi(x)\phi(y)$ and $\phi(x) \geq 1$ when $x \neq 0$,  for $x,y$ in $A$. If 
, for some $s$ in $[0,1]$, we have $\phi(\alpha_i) \geq \phi(\gamma)^{s}$  for all $i$, then the question is to obtain a lower bound for $\sup_{1\leq i < j \leq m}\phi(\alpha_i - \alpha_j)$ in terms of $\phi(\gamma)$, $m$ and $s$. 
This question is relevant, for example, to the problem of determining upper bounds for the number of integer points on small arcs of conics considered in \cite{cil1}, \cite{cil3}, \cite{cil2}, \cite{cil6} and problem of showing that the number of divisors of an integer $N$ lying in certain arithmetical progressions is bounded independently of $N$, considered in \cite{lens}. 

\vspace{2mm}
\noindent
In most situations where the aforementioned question is of interest ({\em loc.\,cit.}), the integral domain $A$ is either a factorial ring or a Dedekind domain and, indeed, it is by assuming that $A$ has one of these properties that this question has been studied.  For instance, when $A$ is a factorial ring we have $\phi(\alpha_i - \alpha_j) \geq \phi((\alpha_i, \alpha_j))$ for $1\leq i < j \leq m$, where $(\alpha_i, \alpha_j)$ is the greatest common divisor of $\alpha_i$ and $\alpha_j$ in $A$. Further, 
(see \cite{cil5}, pages 6 to 8 and also \cite{lens}), there is a natural  measure $\mu$ on the set $X$ of powers of irreducible elements of $A$ dividing $\gamma$ such that , for all distinct $i$ and $j$, we have $\log \phi((\alpha_i , \alpha_j))/\log{\phi(\gamma)} = \mu(E_i \cap E_j)$,  where the $E_i$ are subsets of $X$.  Applying the case of the overlapping theorem of \cite{cil5} that gives a lower bound for $\sum_{1\leq i < j \leq m} \mu(E_i \cap E_j)$ in terms of $\sum_{1\leq i < j \leq m} \mu(E_i)$, one deduces a lower bound for  $\sup_{1\leq i < j \leq m}\phi((\alpha_i, \alpha_j))$, and {\em a forteriori} for $\sup_{1\leq i < j \leq m}\phi(\alpha_i - \alpha_j)$, in terms of $\phi(\gamma), m$ and $s$.

\vspace{2mm}
\noindent
When $A$ is a Dedekind domain, a closely related argument is provided in Theorems 1.1 and 1.2 of \cite{cil2}, based on the observation that the ideal $<\alpha_i -\alpha_j>$ is contained in the ideal $({\mathfrak a}_i, {\mathfrak a}_j)$, which is the greatest common divisor of the ideals  ${\mathfrak a}_i$ and  ${\mathfrak a}_j$, the ideals of $A$ generated respectively by $\alpha_i$ and $\alpha_j$, and assuming that $\phi$ has a natural extension to the ideals of $A$.

\vspace{2mm}
\noindent
In Section 2 we present a simple identity for the Vandermonde determinant that immediately yields, for any integral domain $A$, a lower bound for $\sup_{1\leq i < j \leq m}\phi(\alpha_i - \alpha_j)$ in terms of $\phi(\gamma)$, without recourse to factorization in $A$. This lower bound allows us to easily recover a number of results given in \cite{cil2} and \cite{cil5}. In Section 3 we show  that the case of the overlapping theorem of \cite{cil5} that gives a lower  bound for  $\sup_{1\leq i < j \leq m}\phi((\alpha_i, \alpha_j))$ when $A$ is a factorial ring and Theorem 1.1 of \cite{cil2}, which gives the analogous result when $A$ is a Dedekind domain, may also, in principle, be deduced from the identity given here. We conclude with some notes related to the contents of this article in Section 4.

\end{nsc}

\begin{nsc}{An Identity for the Vandermonde Determinant}

\noindent
Throughout this article $m$ shall denote an integer $\geq 2$.

\vspace{2mm}
\vspace{2mm}
\begin{th}
Let $A$ be a commutative ring and $\{\alpha_i\}_{1 \leq i \leq m}$ and $\{\beta_i\}_{1 \leq i \leq m}$ be sequences of $m$ elements in $A$ for which there is exists a $\gamma$ in $A$ satisfying  $\alpha_i \beta_i = \gamma$ for all $i$. For each integer $k$ satisfying $0 \leq k \leq m-1$ we then have

\vspace{-3mm}
\begin{equation}
\label{det}
 \gamma^{\frac{k(k+1)}{2}}\,\prod_{1 \leq i < j \leq m} (\alpha_i - \alpha_j)
 =
\prod_{1 \leq i \leq m} \alpha_{i}^{k}  \;
\left| \begin{matrix}
{\beta_1}^{k} & {\beta_2}^k & \dots & {\beta_{m}}^k\\
\vdots & \vdots & & \vdots \\
{\beta_1} & {\beta_2} & \dots & {\beta_{m}}\\
1 & 1 & \dots & 1\\
{\alpha_1} & {\alpha_2} & \dots & \alpha_{m}\\
\vdots & \vdots & & \vdots \\
{\alpha_1}^{m-k-1} & {\alpha_2}^{m-k-1} & \dots & {\alpha_{m}}^{m-k-1}
\end{matrix} \right | \; .
\end{equation}

\end{th}

\vspace{4mm}
\noindent
{\sc Proof. ---} When $1 \leq k \leq m-1$ and for each $i$, $1 \leq i \leq m$, we multiply the $i$ th column of the determinant on the right hand side of (\ref{det}) by $\alpha_{i}^k$ . For $1 \leq i \leq m$ and $1\leq j \leq k$ the $(i,j)$ th entry in the resulting determinant is $\beta_i^{k-j+1}\alpha_{i}^{k}= (\beta_i \alpha_i)^{k-j+1}\alpha_{i}^{j-1} = \gamma^{k-j+1}\alpha_i^{j-1}$. Therefore $\gamma^{k-j+1}$ is common to each entry in the $j$th row, for $1\leq j \leq k$. Since 
$\prod_{1 \leq j \leq k} \gamma^{k-j+1} = \gamma^{\frac{k(k+1)}{2}}$, (\ref{det}) now follows on   using the well known evaluation of the Vandermonde determinant, to which it reduces when $k=0$.

\vspace{3mm}
\noindent
{\sc Definition 1. ---} When $A$ is a commutative ring and $\{\alpha_i\}_{1 \leq i \leq m}$ and $\{\beta_i\}_{1 \leq i \leq m}$ are sequences of elements of $A$ we write ${\rm det}_{k}(\alpha,\beta)$, for each integer $k$ satisfying $0\leq k \leq m-1$, to denote the determinant on the right hand side of (\ref{det}).

\vspace{2mm}
\noindent
The preceding definition allows us to  rewrite the identity (\ref{det}) in the following form. For all integers $k$ satisfying $0 \leq k \leq m-1$ and  $\{\alpha_i\}_{1 \leq i \leq m}$, $\{\beta_i\}_{1 \leq i \leq m}$ and $\gamma$ as in Theorem 1 we have 

\vspace{-3mm}
\begin{equation}
\label{det1}
 \gamma^{\frac{k(k+1)}{2}}\,\prod_{1 \leq i < j \leq m} (\alpha_i - \alpha_j)
 \;=
\;{\rm det}_{k}(\alpha, \beta) \;
\prod_{1 \leq i \leq m} \alpha_{i}^{k}  \; .
\end{equation}

\vspace{3mm}
\noindent
In order to choose optimal values of $k$ in applications of (\ref{det1}), we define
, for any real number $s$ in $[0,m]$,

\vspace{-4mm}
\begin{equation}
\label{Ksm}
K(s,m) = \sup_{0 \leq k \leq m-1} (sk -\frac{k(k+1)}{2}) \; .
\end{equation}

\noindent
In this article $K(s,m)$ plays essentially the \/same role as $E_{k}(\gamma)\binom{k}{2}$ in \cite{cil2}, Theorem 1.1 and, by $(i)$ of Lemma~1 below, the same role as $Q_{2}(x)$ in \cite{cil5}.

\vspace{3mm}
\begin{le}
We have the following relations for $K(s,m)$.

\vspace{2mm}
\noindent
(i) For all  $s$ in $[0,m]$ we have $K(s,m) \; = \; \left(s[s] - \frac{[s]([s]+1)}{2}\right) \; \geq \; \frac{s(s-1)}{2} \;$.

\vspace{2mm}
\noindent
(ii) $\frac{K(m/2, m)}{\binom{m}{2}} \; =\; \frac{1}{4}-\frac{1}{8[\frac{m}{2}]+4}$ when
$m$ is an odd integer.

\vspace{2mm}
\noindent
(iii) When $m$ is an integer $\geq 2$, and for all $s$ in $[0,1]$, we have $\frac{K(sm, m)}{\binom{m}{2}} \, \geq \, s^2 - \frac{s(1-s)}{m-1} \, \geq \, s^2 - \frac{1}{4(m-1)}. $

\end{le}

\vspace{4mm}
\noindent
{\sc Proof. ---} Let us verify $(i)$. The function $f(t) = st -\frac{t(t+1)}{2}= (s-\frac{1}{2})t-\frac{t^2}{2}$ is a smooth strictly concave function on ${\bf R}$ that satisfies $f(s) =f(s-1)$. The supremum of $f(t)$ over the integers in $[0, m-1]$ is therefore attained at an integer in $[0, m-1] \cap [s-1, s]$. When $s$ is not  an integer, $[s]$ is the unique integer in this intersection and the required       supremum is attained at $[s]$. When $s$ is an integer, $s = [s]$ and $s-1$ are the integers in  $[0, m-1] \cap [s-1, s]$ and, since $f(s) = f(s-1)$, we see that the required supremum is attained at $[s]$ in this case as well. Moreover, we also have $f([s]) \geq f(s) = s(s-1)/2$. We set $m = 2k+1$ and $s =m/2$ in the equality in $(i)$ and obtain $K(m/2,m) = k^{2}/2$, from which $(ii)$ follows on dividing by $\binom{m}{2}$ and rearranging terms. We obtain $(iii)$ from the inequality in $(i)$ on noting that $s(1-s) \leq 1/4$ when $s$ is in $[0,1]$.

\vspace{3mm}
\begin{prop}
Let $A$ be  an integral domain and  $\alpha = \{\alpha_i\}_{1 \leq i \leq m}$ and $\beta =  \{\beta_i\}_{1 \leq i \leq m}$ be sequences of distinct non-zero elements of $A$. If $\alpha_i \beta_i = \gamma$ for some $\gamma$ in $A$ and for each $i$, then ${\rm det}_{k}(\alpha, \beta)$ is a non-zero element of $A$ for all $k$, $0 \leq k \leq m-1$. 
Suppose that $\phi$ is a function from $A$ into ${\bf R}_{+}$ satisfying $\phi(xy) = \phi(x)\phi(y)$ and $\phi(x) \geq 1$ when $x \neq 0$, for all $x,y$ in $A$. Then 
 $\phi({\rm det}_{k}(\alpha, \beta)) \geq 1$ for all $k$, $0 \leq k \leq m-1$. 

\vspace{2mm}
\noindent
Suppose that $L \geq 0$ satisfies $\phi({\rm det}_{k}(\alpha, \beta)) \geq L$, for all $k$, $0 \leq k \leq m-1$, and that for some $s$ in $[0,1]$ we have $\phi(\alpha_i) \geq \phi(\gamma)^{s} $ for all $i$. We then have 

\vspace{-5mm}
\begin{equation}
\label{ineq}
\sup_{1 \leq i < j \leq m} \phi(\alpha_i- \alpha_j) \; \geq \; L^{\frac{1}{\binom{m}{2}}} \phi(\gamma)^{\frac{K(s,m)}{\binom{m}{2}}} \;. 
\end{equation}

\end{prop}

\vspace{3mm}
\noindent
{\sc Proof. ---} Since $A$ is an integral domain and  $\alpha$,  $\beta$  are sequences of distinct non-zero elements of $A$, we have $\gamma \neq 0$. The left hand side of (\ref{det1}) is thus distinct from 0 and therefore  ${\rm det}_{k}(\alpha, \beta)$ is distinct from 0 for all $k$, $0 \leq k \leq m-1$. Consequently, $\phi({\rm det}_{k}(\alpha, \beta)) \geq 1$ for all $k$, $0 \leq k \leq m-1$. To verify (\ref{ineq}) we apply $\phi$ to both sides of (\ref{det1}) and obtain

\vspace{-3mm}
\begin{equation}
\label{ineq1}
\phi(\gamma)^{\frac{k(k+1)}{2}}\,(\sup_{1 \leq i < j \leq m} \phi(\alpha_i- \alpha_j))^{\binom{m}{2}} \; \geq \phi(\gamma)^{\frac{k(k+1)}{2}}\,\prod_{1 \leq i < j \leq m} \phi(\alpha_i- \alpha_j) \; \geq \; L \, \phi(\gamma)^{sk} \; ,
\end{equation}

\noindent
for all integers $k$, $0 \leq k \leq m-1$.  On rearranging the terms in  (\ref{ineq1}) and using (\ref{Ksm}) we obtain (\ref{ineq}).

\vspace{2mm}
\noindent
The following corollary to Proposition 1 is implicit in \cite{cil2}, proof of Theorem 1.2, where only the case of this corollary for quadratic extensions of ${\bf Q}$ is required and this is obtained in  \cite{cil2} by an application of Theorem 1.1 of \cite{cil2}.  

\vspace{3mm}
\begin{cor}
Suppose that $K$ is number field of degree $n$ over ${\bf Q}$ and that $\{\alpha_i\}_{1 \leq i \leq m}$ is a sequence of distinct non-zero elements of the ring of integers $A$ of $K$. Let ${\mathcal N}(x)$ denote the norm of an element $x$ of $K$. If for each $i$ we have $|{\mathcal N}(\alpha_i)| = R$ then  

\vspace{-4mm}
\begin{equation}
\label{nf}
\sup_{1\leq i < j \leq m}|\,{\mathcal N}(\alpha_i - \alpha_j)\,|^{\frac{1}{n}} \; \geq \; R^{\frac{K(\frac{m}{n},m)}{\binom{m}{2}}} \; .
\end{equation}

\end{cor}

\vspace{2mm}
\noindent
{\sc Proof. ---} Since $|{\mathcal N}(\alpha_i)| = R$ for each $i$, $R$ belongs to the ideal generated by each $\alpha_i$ in $A$. Thus on setting $\gamma = R$, there exists, for each $i$, a $\beta_i$ in $A$ such that $\alpha_i \beta_i = \gamma$. Let  $\phi$ be the function $x \rightarrow  |{\mathcal N}(x)|^{\frac{1}{n}}$. Since $R$ is in ${\bf Z}$, we have $\phi(R) = R$ and hence $\phi(\alpha_i) = \phi(\gamma)^{\frac{1}{n}}$ for all $i$. The corollary now follows from Proposition 1 applied with $L = 1$ and $s = 1/n$.

\vspace{2mm}
\noindent
The following corollary to Proposition 1 is implicit in the proof of Proposition 1 of H. Lenstra \cite{lens}, whose methods are closely related to the case of the overlapping theorem of \cite{cil5} mentioned in Section 1 above.

\vspace{2mm}
\begin{cor}
Let $s$ be a real number in $(0,1)$ and $\{d_i\}_{1 \leq i \leq m}$ be distinct positive divisors of an integer $N\geq 1$ and satisfying $d_i \geq N^{s}$ for all $i$. If each $d_i$ belongs to the arithmetic progression $a \,{\rm mod}\, q$, where $(a,q)=1 $, we then have  

\vspace{-4mm}
\begin{equation}
\label{lens1}
\sup_{1 \leq i  < j \leq m} |d_i - d_j| \,\geq \, q\, N^{\frac{K(sm,m)}{\binom{m}{2}}} \; .
\end{equation}
\end{cor}

\vspace{2mm}
\noindent
{\sc Proof. ---} We take $A = {\bf Z}$ and set $\alpha_i = d_i$, $\beta_i = \frac{N}{d_i}$ and $\gamma = N$ and take $\phi$ to be the function $ x \rightarrow |x|$. Since each $\alpha_i \equiv a \,{\rm mod}\, q$, we see that $\prod_{1 \leq i < j \leq m} (\alpha_i -\alpha_j)$ is divisible by $q^{\binom{m}{2}}$. Since $(a, q) = 1 $, we see that $\prod_{1 \leq i \leq m} \alpha_i^{k} \not \equiv 0 \, {\rm mod}\, q$, for any integer $k \geq 0$. The identity (\ref{det1}) then shows that ${\rm det}_k(\alpha, \beta)$ is divisible by $q^{\binom{m}{2}}$, for all integers $k$, $0 \leq k \leq m-1$, and hence that we may take $L = q^{\binom{m}{2}}$ when applying Proposition 1.

\vspace{2mm}
\noindent
The following corollary to Proposition 1 generalises Theorem 1.4 of \cite{cil2}. 

\vspace{2mm}
\begin{cor}
Suppose that $E$ is an integral domain and $X = (X_{\iota})_{{\iota} \in I}$ is a family of indeterminates indexed by a set $I$. Let $\{P_i(X)\}_{1 \leq i \leq m}$ be a sequence of distinct polynomials in $E[X]$. If $R(X)$ is a common multiple of the polynomials $P_i(X)$ in $E[X]$ and if, for some $s$ in $[0,1]$,  ${\rm deg}(P_{i}) \geq s\, {\rm deg}(R)$ for all $i$, we then have 

\vspace{-4mm}
\begin{equation}
\label{lens11}
\sup_{1 \leq i  < j \leq m} {\rm deg}(P_{i} - P_{j}) \;\geq \; {\rm deg}(R)\,\frac{K(sm,m)}{\binom{m}{2}} \;, 
\end{equation}

\noindent
where ${\rm deg}(u)$ denotes the total degree of a polynomial $u(X)$ in $E[X]$.

\end{cor}
 
\vspace{2mm}
\noindent
{\sc Proof. ---} Since $E$ is an integral domain so is $E[X]$ and ${\rm deg}(uv) = {\rm deg}(u) + {\rm deg}(v)$ for $u$ and $v$ elements of $E[X]$. We apply Proposition 1 with  $A = E[X]$, $\alpha_i = P_{i}(X)$,  $\beta_i = Q_{i}(X)$ such that $P_{i}(X)Q_{i}(X) = R(X)$,  $\gamma = R(X)$,  $\phi$ taken to be the function $u \rightarrow  {\rm exp}({\rm deg}(u))$ and $L =1$.

\vspace{2mm}
\noindent
Corollary 1 is the essential point in the proof of Theorem 1.2 of \cite{cil2}, which contains Theorem~1 of \cite{cil1} and improves on the main results of \cite{cil3}, \cite{cil6}. Corollary 2 is the essential point in the proof of Proposition 1 of \cite{lens} as well as Lemma 3.1 of \cite{regis}. We restrict ourselves here to giving only a proof of a version of Theorem 1.2 of \cite{cil2} refering the reader to pages 6 to 8 of \cite{cil5} for an account of the other results.

\vspace{2mm}
\begin{th}
When $d \neq 0, -1$ is a squarefree integer and $m$, $R$ are integers with $m \geq 2$, there are no more than $m$ integer points on any arc of length $\leq\,\frac{|R|^{s(m)}}{\sqrt{|d|}}$ on the conic 

\vspace{-3mm}
\begin{equation}
\label{con}
X^2 + dY^2 = R,
\end{equation}

\noindent
where $s(m) =  1/4 - 1/(8[m/2]+4)$.
\end{th}

\vspace{3mm}
\noindent
{\sc Proof. ---} Indeed, if $\{p_i\}_{1 \leq i \leq m}$ is a sequence of $m$ integer points $p_i = (x_i, y_i)$ on (\ref{con}) and, for each $i$, $\alpha_i = x_i + \sqrt{-d}\,y_i$, then $\alpha_i$ are elements of the ring of integers of ${\bf Q}(\sqrt{-d})$. Since $d$ is a squarefree integer $ \neq 0, -1$, ${\bf Q}(\sqrt{-d})$ is a quadratic extension of ${\bf Q}$ and the triangle inequality gives 

\vspace{-1mm}
\noindent
\begin{equation}
\label{norms}
|d|\; {\| p_i -p_j \|^{2}_{2}} \; \geq \; |{\mathcal N}(\alpha_i -\alpha_j)| \; ,
\end{equation}

\noindent
for all $(i,j)$, where $\|\;\|_2$ denotes the Euclidean distance and ${\mathcal N}$ the norm on  ${\bf Q}(\sqrt{-d})$. If the points $\{p_i\}_{1 \leq i \leq m}$ lie on an arc of length $l$, we have $l  > \| p_i -p_j \|_{2}$ for all $(i,j)$. Since ${\mathcal N}(\alpha_i) = R $ for each $i$, it then follows from (\ref{norms}) and Corollary 1 applied with $n =2$ that

\vspace{-1mm}
\noindent
\begin{equation}
\label{i}
|d|^{\frac{1}{2}}\; l \; > \; \sup_{1 \leq i < j \leq m} |{\mathcal N}(\alpha_i -\alpha_j)|^{\frac{1}{2}} \;  \geq \; |R|^{\frac{K(\frac{m}{2},m)}{\binom{m}{2}}} \; = \; |R|^{s(m)},
\end{equation}

\noindent
when $m$ is an odd integer $\geq 2$, where the equality follows from $(ii)$ of Lemma 1. Plainly, (\ref{i}) implies that there are no more than $m-1$ integer points on an arc of length  $\, \leq \,\frac{|R|^{s(m)}}{\sqrt{|d|}}$ when   $m$ is an odd integer $\geq 2$. When $m$ is an even integer $\geq 2$ we note that $s(m) = s(m+1)$ and apply the preceding conclusion to $m+1$.

\vspace{3mm}
\noindent
{\sc Remark 1. ---} Theorem 1.2 in \cite{cil2} states that when $d \neq 0,1$ is a fixed squarefree integer, on the conic  $X^2 -dY^2 = N$, an arc of length  $N^{\alpha}$ with $\alpha = 1/4 - 1/(8[k/2]+4)$ contains at most $k$ lattice points. This statement appears to be inaccurate with regard to the dependence of the lengths of the arcs on $d$. As Example 1 below shows, there are infinitely many integers $R \geq 1$ such that there are arcs of length  $\frac{2^{13/6} R^{1/6}}{d^{1/3}}$ containing 3 integer points on the ellipses $X^2 + d Y^2 = R^2$, for any integer $d \geq 1$, while Theorem~1.2 of  \cite{cil2} implies that there are no more than $2$ integer points on any arc of length $R^{1/6}$ on these conics.  

\vspace{3mm}
\noindent
{\sc Remark 2. ---} The dependence of the lengths of the arcs on $d$ given by Theorem 2 may be improved  by noting that $x^2 \equiv \, R \,{\rm{mod}}\, p$, for all integer points $(x,y)$ on $X^2 + dY^2 = R$ and primes $p$ dividing $d$. We explain this using the notation of the proof of Theorem 2. Let us first verify that for any prime $p$ dividing $d$ we have $v_p( \prod_{1 \leq i < j \leq m} {\mathcal N}(\alpha_i -\alpha_j)) \geq [\frac{1}{2}(\frac{m^2}{2}-m)]+1$, which we denote by $t(m)$. Indeed, if $k$ of the $x_i$ belong to the same residue class modulo $p$, we then have that $v_p( \prod_{1 \leq i < j \leq m} {\mathcal N}(\alpha_i -\alpha_j)) \geq k(k-1)/2$. Since $x_i^2 \equiv \, R \,{\rm{mod}}\, p$, each $x_i$ lies in one of no more than  2 residue classes modulo $p$. It then follows that for some integer $k$, $0 \leq k \leq m$, we have   

\vspace{-3mm}
\begin{equation}
\label{v_p}
v_p( \prod_{1 \leq i < j \leq m} {\mathcal N}(\alpha_i -\alpha_j)) \; \geq \; \frac{k(k-1)}{2}\, + \, \frac{(m-k)(m-k-1)}{2} \; \geq \; t(m) \, .
\end{equation}

\noindent
Suppose that $p$ divides $d$ but not $R$. Then the identity (\ref{det1}) shows that $v_p({\mathcal N}({\rm det}_{k}(\alpha, \beta)))$ is the same as  $v_p( \prod_{1 \leq i < j \leq m} {\mathcal N}(\alpha_i -\alpha_j))$ and therefore  $v_p({\mathcal N}({\rm det}_{k}(\alpha, \beta))) \geq t(m)$, for such primes $p$. This  bound may be seen to be valid even when $p$ divides $d$ and $R$. In effect, in this case each of the ideals $<\alpha_i>$ and $<\beta_i>$, generated in the ring $A$ of integers of ${\bf Q}(\sqrt{-d})$ by $\alpha_i$ and $\beta_i$ respectively, is divisible by the prime ideal ${\mathfrak p}$, the unique prime ideal lying above the ramified prime $p$ in ${\bf Q}(\sqrt{-d})$. On expanding the determinants ${\rm det}_k(\alpha, \beta)$ with respect to any row, we see that for all integers $k$, $0 \leq k \leq m-1$, we have 

\vspace{-3mm}
\begin{equation}
\label{vp2}
v_{\mathfrak{p}}(<{\rm det}_{k}(\alpha, \beta)>) \; \geq \; \frac{k(k+1)}{2}\, + \, \frac{(m-1-k)(m-k)}{2} \; \geq \; t(m) \, , 
\end{equation}

\noindent
where $<{\rm det}_{k}(\alpha, \beta)>$ is the ideal generated by ${\rm det}_{k}(\alpha, \beta)$ in $A$. Consequently, we have $v_p({\mathcal N}({\rm det}_{k}(\alpha, \beta))) \geq t(m)$ even when $p$ divides $d$ and $R$. Since $d$ is a squarefree integer, we then deduce that $|{\mathcal N}({\rm det}_{k}(\alpha, \beta)))| \geq |d|^{t(m)}$. 

\vspace{2mm}
\noindent
On using the bound $|{\mathcal N}({\rm det}_{k}(\alpha, \beta)))| \geq |d|^{t(m)}$ in the proof of Corollary 1 and arguing as in the proof of Theorem 2, we see that $\frac{R^{s(m)}}{\sqrt{|d|}}$ in the statement of Theorem 2 maybe replaced by $ \frac{R^{s(m)}}{|d|^{l(m)}}$, where $l(m)$ is defined to be $ \frac{1}{2}(1- \frac{t(m)}{\binom{m}{2}})$ when $m$ is odd and $l(m) = l(m+1)$ when $m$ is even. In particular, we see that there are no more than 2 integer points on an arc of length $\frac{R^{1/6}}{|d|^{1/3}}$ on the conic $X^2 +dY^2 =R$, with $d$ and $R$ as in Theorem 2.  

\vspace{3mm}
\noindent
The following example was kindly supplied to the author by Prof. Joseph Oesterl{\'e}.

\vspace{2mm}
\noindent
{\sc Example 1. ---} Let $t$ and $d$ be integers $\geq 1$ and let $u = d^2t + dt -d +1$. Let $p_i = (x_i,y_i)$, $1 \leq i \leq 3$, be points in the plane with coordinates $x_i$, $y_i$ given below.

\vspace{-3mm}
\begin{equation}
\label{coord}
\begin{array}{ll}
x_1 = dt(2dt-1)u -1, & y_1 = t(2dt+1)u +1\\
x_2 = x_1 +2dt+2, & y_2 = y_1 -2dt\\
x_3 = x_1 -2dt, & y_3 = y_1 +2dt-2\\
\end{array} 
\end{equation}

\noindent
We then verify that $x_i^2 + dy_i^2 = x_1^2 + dy_1^2$, for $1 \leq i \leq 3$ and, on setting $R = x_1^2 + dy_1^2$, we see that the points $p_i$ are integer points all of which lie on the ellipse $X^2 + d Y^2 = R$. Let us set $D = \sup_{1 \leq i < j \leq 3} \|p_i - p_j\|_2$ and $l$ to be the length of the shortest arc on the ellipse containing the points $p_i$. Then as $t\rightarrow +\infty$ we have 

\vspace{-3mm}
\begin{equation}
\label{ex1}
R \; \sim \; 4d^7(d+1)t^6, \; \;  D \; \sim \; 4\sqrt{2}dt \;\; \text{and} \;\; l \sim D , 
\end{equation}

\noindent
where the relation $l \sim D$ follows on noting that $\frac{D}{R^{\frac{1}{2}}} \rightarrow 0$ as $t \rightarrow +\infty$. Since $d \geq 1$, it follows from (\ref{ex1}) that

\vspace{-4mm}
\begin{equation}
\label{ex2}
l < \frac{2^{13/6}  R^{\frac{1}{6}}}{d^{\frac{1}{3}}} \; \; \text{for all sufficiently large $t$.}
\end{equation}

\vspace{3mm}
\noindent
{\sc Remark 3. ---} On setting $\alpha_i = x_i + \sqrt{-d}y_i$ and $\beta_i = x_i -\sqrt{-d}y_i$ for $1 \leq i \leq 3$, with $x_i$ and $y_i$ as in Example 1,  we see that ${\rm det}_{1}(\alpha, \beta) = 16\sqrt{-d}$, so that the lower bound used in the proof of Corollary 1 for $|{\mathcal N}({\rm det}_k(\alpha,\beta))|$ is best possible with respect to $R$ when $m =3$ and $K = {\bf Q}(\sqrt{-d})$, $d$ an integer $\geq 1$. The author does not know if this lower bound, and, similarly, the lower bounds for $\phi({\rm det}_k(\alpha, \beta))$ used in the  proofs of corollaries 2 and 3 above, may be improved upon for large values of $m$.

\vspace{3mm}
\noindent
It will interest the reader to note that a recent conjecture (Conjecture 14 on page 11 of \cite{cil4}) of J.~ Cilleruelo and A. Granville looks forward to a considerable improvement of Theorem 2 when the conic in this theorem is a circle. On page 15 of \cite{cil4}, Cilleruelo and Granville give a flowchart relating their conjecture to a number of other interesting conjectures on the interface between Fourier Analysis and Number Theory. Also, on page 12 of the same article the reader will find a summary of what is known on the theme of Theorem 2.

\end{nsc}

\begin{nsc}{The Overlapping Theorem, Divisors in a Dedekind Domain and the Identity.}

\begin{prop}
Let $\{a_i\}_{1 \leq i \leq m}$ be a sequence of real numbers with each $a_i \geq 0$. For each integer $k$, $0 \leq k \leq m-1$, we then have the inequality 

\vspace{-4mm}
\begin{equation}
\label{over1}
\frac{k(k+1)}{2}\sup_{1 \leq i \leq m} a_i + \sum_{1 \leq i < j \leq m} \inf(a_i, a_j) \; \geq \; k\sum_{1 \leq i \leq m} a_i \; .
\end{equation}

\end{prop}

\vspace{2mm}
\noindent
{\sc Proof. ---} Suppose first that the $a_i$ are {\em distinct integers} $\geq 0$. Let $p$ be a prime number and let us apply the identity (\ref{det1}) of Section 2, with 
$\alpha_i = p^{a_i}$, $\beta_i= p^{(\sup{a_i}) - a_i}$ and $\gamma = p^{\sup{a_i}}$ for $1 \leq i \leq m$. Then ${\rm det}_{k}(\alpha, \beta)$ is an integer distinct from 0 for all $k$, $0 \leq k \leq m-1$. We now obtain (\ref{over1}) on comparing the powers of $p$ dividing both side of (\ref{det1}) of Section 2 and noting that, since the $a_i$ are distinct,  $v_{p}(p^{a_{i}} - p^{a_{j}}) = \inf(a_i,a_j)$ for all $(i,j)$, $1 \leq i < j \leq m$. 

\vspace{2mm}
\noindent
When the $a_i$ are {\em distinct rational numbers} $\geq 0$, we write them to a common denominator, apply (\ref{over1}) to their numerators, which are then distinct integers $\geq 0$, and divide by throughout by their common denominator. Finally, noting that the set of points $(a_1, a_2, \ldots, a_m)$ in  ${\bf R}^m$ with $a_i$ distinct rational numbers $\geq 0$ is dense in the subset of ${\bf R}^m$ consisting of $(a
_1, a_2, \ldots, a_m)$, with each $a_i \geq 0$, we obtain (\ref{over1}) by continuity.

\vspace{2mm}
\noindent
{\sc Remark 1. ---} The inequality (\ref{over1}) may evidently be verified directly as well by reducing to the case when the $a_i$ are in increasing order and comparing the two sides as in the proof of Theorem 1.1 in \cite{cil2}.

\vspace{5mm}
\noindent
The case of the overlapping theorem of \cite{cil5} mentioned  in Section 1 is the following corollary to Proposition 1, which is stated using the notation $K(s,m)$ of (\ref{Ksm}) of Section 2. From a conceptual point of view, the proof of the following corollary is closely related to that in \cite{cil5}.

\begin{cor}
When $X$ is a  measure space with a probability measure $\mu$ and $\{A_i\}_{1 \leq i \leq m}$ is a finite sequence of measurable subsets of $X$ we have the inequality

\vspace{-4mm}
\begin{equation}
\label{over2}
 \sum_{1 \leq i < j \leq m} \mu(A_i \cap A_j) \; \geq \;  K\left(\sum_{1 \leq i \leq m} \mu(A_i), m\right) \; .
\end{equation}

\end{cor}

\vspace{2mm}
\noindent
{\sc Proof. ---} For each $t$ in $X$ we apply (\ref{over1}) to $\{\chi_{i}(t)\}_{1 \leq i \leq m}$, where the $\chi_{i}$ are the characteristic functions of the sets $A_i$. On noting that $t \rightarrow \sup_i{\chi_{i}}(t)$ is the characteristic function of $\cup_{1 \leq i \leq m}A_i$ and that $t \rightarrow \inf(\chi_{i}(t), \chi_{j}(t))$ is the characteristic function of $A_i \cap A_j$ and integrating the resulting relation with respect to $\mu$ we obtain, for every integer $k$ satisfying $0 \leq k \leq m-1$, that

\vspace{-5mm}
\begin{equation}
\label{over4}
\frac{k(k+1)}{2} \mu(\cup_{1\leq i\leq m} A_i) + \sum_{1 \leq i < j \leq m} \mu(A_i \cap A_j) \; \geq \; k\sum_{1 \leq i \leq m} \mu(A_i) \; .
\end{equation}

\noindent
Since $ \mu(\cup_{1\leq i\leq m} A_i) \leq 1$ and $0 \leq \sum_{1 \leq i \leq m} \mu(A_i) \leq m$, we conclude using (\ref{Ksm}) of Section 1.  

\vspace{4mm}
\noindent
The following corollary to Proposition 1 is Theorem 1.1 of \cite{cil2}. For the sake of completeness we give a proof, which is the same as given in  \cite{cil2}. 

\vspace{3mm}
\begin{cor}
Suppose that $A$ is a Dedekind domain and that $\{{\mathfrak a}_i\}_{1 \leq i \leq m}$ is a sequence of non-zero ideals in $A$. Suppose that ${\mathfrak c}$ is a non-zero ideal in $A$ which is divisible by each of the ${\mathfrak a}_i$ then for each integer $k$ satisfying $0 \leq k \leq m-1$ we have  

\vspace{-4mm}
\begin{equation}
\label{over3}
{ {\mathfrak c}}^{\frac{k(k+1)}{2}} \prod_{1 \leq i < j \leq m} ({\mathfrak a}_i, {\mathfrak a}_j) \;\; \text{is divisible by} \;\; \prod_{1 \leq i \leq m} {{\mathfrak a}_{i}}^k \; , 
\end{equation}

\noindent
where $(\,,\,)$ denote the greatest common divisor. Consequently, when $\phi$ is a function from the ideals set of $A$ into ${\bf R}_{+}$ satisfying $\phi(\mathfrak{a}\mathfrak{b}) = \phi(\mathfrak{a})\phi(\mathfrak{b})$ and $\phi(\mathfrak{a}) \geq 1$ when $\mathfrak{a} \neq 0$,  for ideals $\mathfrak{a}, \mathfrak{b}$ in $A$ and, if
, for some $s$ in $[0,1]$, we have $\phi(\mathfrak{a}_i) \geq \phi(\mathfrak{c})^{s}$  for all $i$, then 

\vspace{-3mm}
\begin{equation}
\label{t1.1}
\sup_{1 \leq i < j \leq m} \phi(({\mathfrak a}_i, {\mathfrak a}_j)) \; \geq \; \phi(\mathfrak{c})^{\frac{K(sm,m)}{\binom{m}{2}}} \; .
\end{equation}

\end{cor}

\vspace{2mm}
\noindent
{\sc Proof. ---} Since the ${\mathfrak a}_i$ are ideals in $A$, we have $v_{{\mathfrak p}}({\mathfrak a}_i) \geq 0$ for all prime ideals ${\mathfrak p}$ in $A$ and all $i$. Since each ${\mathfrak a}_i$ divides  ${\mathfrak b}$ we have $ v_{\mathfrak p}({\mathfrak c}) \geq \sup_{i} v_{\mathfrak p}({\mathfrak a}_i)$ for all prime ideals ${\mathfrak p}$ in $A$ and all $i$. On comparing the exponents of ${\mathfrak p}$ in the two expressions in (\ref{over3}) we then see that (\ref{over3}) follows from (\ref{over1}) applied to $\{ v_{\mathfrak p}({\mathfrak a}_i) \}_{1 \leq i \leq m }$,  for each prime ideal ${\mathfrak p}$ in $A$. The properties of $\phi$ and  (\ref{over3}) imply that 

\vspace{-3mm}
\begin{equation}
\label{t1.11}
\phi(\mathfrak{c})^{\frac{k(k+1)}{2}}\,\prod_{1 \leq i < j \leq m} \phi(({\mathfrak a}_i, {\mathfrak a}_j)) \; \geq \; \prod_{1 \leq i \leq m}\phi(\mathfrak{a}_i)^{k} \; ,
\end{equation}

\noindent
for every integer $k$, $0 \leq k \leq m-1$, from which (\ref{t1.1}) follows in the manner of the proof of Proposition 1 of Section 2. 

\vspace{3mm}
\noindent
{\sc Remark 2. ---} When $A$ is a principal ideal domain,  (\ref{t1.1}) shows that the lower bound for $\sup_{1 \leq i < j \leq m} \phi(a_i - a_j)$ provided by Proposition 1 of Section 1 applied with $L =1$ is, in fact, a lower bound for $\sup_{1 \leq i < j \leq m} \phi((a_i,a_j))$. This conclusion may be obtained for any factorial ring $A$ by using (\ref{over2}) in place of (\ref{t1.1}), as described in Section 1. For each integer $m \geq 2$, there are examples that show the inequalities (\ref{over2}) and (\ref{t1.1}) cannot be improved in general (see Theorems 2.2 and 3.7 of \cite{cil5}).

\end{nsc}

\begin{nsc}{notes}

\vspace{2mm}
\noindent
The author arrived at the identity (\ref{det}) of Section 1 as one way of generalising the elementary formula $abc = 4\Delta R$, where $a$,$b$ and $c$ are the sides of a triangle, $\Delta$ its area and $R$, the radius of its circumcircle. Indeed, if one applies the identity with $m =3$, $k =1$, $\alpha_i$ elements of ${\bf C}$ denoting the vertices of the triangle, $\beta_i = \bar{\alpha_i}$, $\gamma = R^2$, one arrives at the formula  $abc = 4\Delta R$ on taking absolute values of both sides of the resulting relation and noting that $|{\rm det}_1(\alpha, \beta)| = 4\Delta$. The use of the formula $abc = 4\Delta R$ in obtaining the case of Theorem 2 of Section 2 when $m = 2$ and when the conic in this theorem is a circle is described on page 899 of \cite{cil1}.

\vspace{2mm}
\noindent
The use of a relation between matrices of the form $(f_i(x_j))$ and $(x_j^{i-1})$, where $x_j$ are elements of a commutative ring $A$ - usually a subring of the complex numbers - and $f_i$ suitable functions on this ring, to study the gaps between the $x_j$ is well known in the context of the Bombieri-Pila method. Indeed, even the simplest of such relations, namely the case when the $f_i$ are polynomials, may be used to deduce interesting conclusions, as for example, in the second proof of Theorem 10 on page 7 of \cite{cil4}; the identity (\ref{det}) of Section 1 may certainly be viewed from this perspective as well.
Also, the reader will not miss the close relation between the method of proof of this identity and K. Mahler's manipulation of the Vandermonde determinant in the proof of his well known upper bound for the discriminant of a polynomial in \cite{mah}.

\vspace{2mm}
\noindent
Finally, we note that there are applications described in \cite{cil5} of even the particular case of the overlapping theorem that we have been concerned with here on which the identity of this article does not shed any light.

\end{nsc}

\vspace{4mm}
\noindent
{\bf Acknowledgment :} Through the course of preparing this article the author benefitted considerably from correspondence with Prof. J. Cilleruelo, from whose works, in particular those with Prof. G. Tenenbaum and Prof. A. Granville, the author learnt of a number of results relevant to the contents here. The author expresses his gratitude to these scholars as well as to Prof. O. Ramar{\'e}, Dr. T.D. Browning and Dr. Gyan Prakash, who very  patiently went through a number of versions of this article and provided the author with several useful suggestions.

\vspace{-2mm}
\bibliographystyle{plain}

\vspace{1cm}
\begin{flushleft}

{\em
Harish-Chandra Research Institute, \\
Chhatnag Road, Jhunsi,\\
Allahabad - 211 019, India.}\\

{\em email} : suri@mri.ernet.in
\end{flushleft}

\end{document}